\documentclass[11pt,reqno,a4paper,twoside]{article}
\usepackage{amsmath,amsthm,amstext,amscd,amssymb,euscript}

\usepackage[english]{babel}
 \makeatletter
 \@addtoreset{equation}{section}
 \makeatother

 \makeatletter
 \@addtoreset{enunciato}{section}
 \makeatother

 \newcounter{enunciato}[section]

 \newtheorem{ittheorem}{Theorem}
 \newtheorem{itlemma}{Lemma}
 \newtheorem{itproposition}{Proposition}
 \newtheorem{itdefinition}{Definition}

 \newenvironment{theorem}{\addtocounter{enunciato}{1}
 \begin{ittheorem}}{\end{ittheorem}}

 \newenvironment{lemma}{\addtocounter{enunciato}{1}
 \begin{itlemma}}{\end{itlemma}}

 \newenvironment{proposition}{\addtocounter{enunciato}{1}
 \begin{itproposition}}{\end{itproposition}}

 \newenvironment{definition}{\addtocounter{enunciato}{1}
 \begin{itdefinition}}{\end{itdefinition}}

 \newcommand{\be}[1]{\begin{equation}\label{#1}}
 \newcommand{\ee}{\end{equation}}

 \newcommand{\bt}[1]{\begin{theorem}\label{#1}}
 \newcommand{\et}{\end{theorem}}

 \newcommand{\bl}[1]{\begin{lemma}\label{#1}}
 \newcommand{\el}{\end{lemma}}

 \newcommand{\bp}[1]{\begin{proposition}\label{#1}}
 \newcommand{\ep}{\end{proposition}}

 \newcommand{\bd}[1]{\begin{definition}\label{#1}}
 \newcommand{\ed}{\end{definition}}

 \newcommand{\bpr}{\begin{proof}}
 \newcommand{\epr}{\end{proof}}

\newcommand{\Zd}{\mathbb Z^d}

\renewcommand{\phi}{\varphi}
\renewcommand{\subset}{\subseteq}

\def\1{ {\mathit{1} \!\!\>\!\! I} }

\newcommand{\Z}{\mathbb Z}
\newcommand{\R}{\mathbb R}
\newcommand{\N}{\mathbb N}

\newcommand{\C}{\mathcal C}

\newcommand{\pee}{\mathbb P}
\newcommand{\E}{\mathbb E}
\newcommand{\T}{\mathbb T}

\newcommand{\ce}{\ensuremath{\mathcal{C}}}

\newcommand{\loc}{\ensuremath{\mathcal{L}}}

\newcommand{\la}{\ensuremath{\Lambda}}
\newcommand{\si}{\ensuremath{\sigma}}
\newcommand{\epsi}{\ensuremath{\epsilon}}

\newcommand{\tend}{\ensuremath{\mathbb T}_N^d}
\newcommand{\norm}{{|\tend|}}
\newcommand{\caD}{\mathcal D}
\newcommand{\caP}{\mathcal P}
\newcommand{\caH}{\mathcal H}
\newcommand{\caM}{\mathcal M}
\newcommand{\caI}{\mathcal I}
\newcommand{\dgam}{\dot{\gamma}}
\newcommand{\cM}{\mathcal{M}}
\newcommand{\cN}{\mathcal{N}}

\begin{document}

\title{
A large-deviation view on\\ 
dynamical Gibbs-non-Gibbs transitions}

\author{
\renewcommand{\thefootnote}{\arabic{footnote}}
A.C.D. van Enter
\footnotemark[1] 
\\
\renewcommand{\thefootnote}{\arabic{footnote}}
R. Fern\'andez
\footnotemark[2]
\\
\renewcommand{\thefootnote}{\arabic{footnote}}
F. den Hollander
\footnotemark[3]
\\
\renewcommand{\thefootnote}{\arabic{footnote}}
F. Redig
\footnotemark[4]
}

\footnotetext[1]{
Johann Bernoulli Institute of Mathematics and Computer Science, University of Groningen,
P.O.\ Box 407, 9700 AK, Groningen, The Netherlands, {\sl A.C.D.van.Enter@rug.nl}  
}
\footnotetext[2]{
Department of Mathematics, Utrecht University, P.O.\ Box 80010, 3508 TA Utrecht, 
The Netherlands, {\sl R.Fernandez1@uu.nl}
}
\footnotetext[3]{
Mathematical Institute, Leiden University, P.O.\ Box 9512, 2300 RA, Leiden, The 
Netherlands, \newline {\sl denholla@math.leidenuniv.nl}
}
\footnotetext[4]{
IMAPP, Radboud University Nijmegen, Heyendaalseweg 135, 6525 AJ Nijmegen, The 
Netherlands, {\sl f.redig@math.ru.nl}
}

\maketitle

\begin{abstract}
We develop a space-time large-deviation point of view on Gibbs-non-Gibbs transitions
in spin systems subject to a stochastic spin-flip dynamics. Using the general 
theory for large deviations of functionals of Markov processes outlined in Feng 
and Kurtz~\cite{FeKu06}, we show that the trajectory under the spin-flip dynamics 
of the empirical measure of the spins in a large block in $\Zd$ satisfies a large 
deviation principle in the limit as the block size tends to infinity. The associated 
rate function can be computed as the action functional of a Lagrangian that is the 
Legendre transform of a certain non-linear generator, playing a role analogous to 
the moment-generating function in the G\"artner-Ellis theorem of large deviation 
theory when this is applied to finite-dimensional Markov processes. This rate function 
is used to define the notion of ``bad empirical measures'', which are the discontinuity 
points of the optimal trajectories (i.e., the trajectories minimizing the rate function) 
given the empirical measure at the end of the trajectory. The dynamical Gibbs-non-Gibbs 
transitions are linked to the occurrence of bad empirical measures: for short times 
no bad empirical measures occur, while for intermediate and large times bad empirical 
measures are possible. A future research program is proposed to classify the various 
possible scenarios behind this crossover, which we refer to as a ``nature-versus-nurture'' 
transition.

{\it Acknowledgment:}  The authors are grateful for extended discussions with 
Christof K\"ulske. Part of this research was supported by the Dutch mathematics
cluster \emph{Nonlinear Dynamics of Natural Systems}. RF is grateful to NWO 
(Netherlands) and CNRS (France) for financial support during his sabbatical leave 
from Rouen University in the academic year 2008--2009, which took place at Groningen 
University, Leiden University and EURANDOM. In the Fall of 2008 he was EURANDOM-chair. 
\end{abstract}


\section{Introduction, main results and research program}
\label{S1}


\subsection{Dynamical Gibbs-non-Gibbs transitions}
\label{S1.1}

Since the discovery of the Griffiths-Pearce-Israel pathologies in renormalization-group
transformations of Gibbs measures, there has been an extensive effort towards
understanding the phenomenon that a simple transformation of a Gibbs measure 
may give rise to a non-Gibbs measure, i.e., a measure for which no reasonable 
Hamiltonian can be defined (see van Enter, Fern\'andez and Sokal~\cite{vEnFeSo93}, 
Fern\'andez~\cite{Fe06}, and the papers in the EURANDOM workshop proceedings 
\cite{MPRF}). From the start, R.L. Dobrushin was interested and involved in 
this development; indeed, Dobrushin and Shlosman~\cite{DoSh97}, \cite{DoSh99} 
proposed a programme of \emph{Gibbsian restoration}, based on the idea that the 
pathological \emph{bad configurations} of a transformed Gibbs measure 
(i.e., the essential points of discontinuity of some of its finite-set, 
e.g. single-site, conditional probabilities) 
are exceptional in the measure-theoretic sense (i.e., they form a set of measure zero). 
This has led to two extended notions of Gibbs measures: \emph{weakly Gibbsian} 
measures and \emph{almost Gibbsian} measures (see Maes, Redig and Van 
Moffaert~\cite{MaVaMoRe99}). Later, several refined notions were proposed, such 
as \emph{intuitively weakly Gibbs} (Van Enter and Verbitskiy~\cite{vEnVe04}) 
and right-continuous conditional probabilities.

In Van Enter, Fern\'andez, den Hollander and Redig~\cite{vEnFedHoRe02}, the 
behavior of a Gibbs measure $\mu$ subject to a high-temperature Glauber spin-flip
dynamics was considered. A guiding example is the case where we start from the 
low-temperature plus-phase of the Ising model, and we run a high-temperature 
dynamics, modeling the fast heating up of a cold system. The question of 
Gibbsianness of the measure $\mu_t$ at 
time $t>0$ can then be interpreted as the existence of a reasonable notion of an
\emph{intermediate-time-dependent temperature} (at time $t=0$ the temperature is 
determined by the choice of the initial Gibbs measure, while at time $t=\infty$ 
the temperature is determined by the unique stationary measure of the dynamics). 
For infinite-temperature dynamics, the effect of the dynamics is simply that of 
a single-site Kadanoff transformation, with a parameter that depends on time. 
The extension to high-temperature dynamics was achieved with the help of a 
\emph{space-time cluster expansion} developed in Maes and Netocn\'y~\cite{MaNe02}. 
The basic picture that emerged from this work was the following:
\begin{itemize}
\item[(1)] 
$\mu_t$ is Gibbs for small $t$;
\item[(2)] 
$\mu_t$ is non-Gibbs for intermediate $t$;
\item[(3)]
in zero magnetic field $\mu_t$ remains non-Gibbs for large $t$, while in non-zero 
magnetic field $\mu_t$ becomes Gibbs again for large $t$.
\end{itemize}

Further research went into several directions and, roughly summarized, gave the 
following results:
\begin{itemize}
\item[(a)] 
Small-time conservation of Gibbsianness is robust: this holds for a large class of spin systems and of dynamics, including discrete spins (Le Ny and Redig~\cite{LNRe02}), continuous 
spins (Dereudre and Roelly~\cite{DeRo05}, van Enter, K\"ulske, Opoku and Ruszel,
~\cite{KuOp08a}, ~\cite{vEnRu08}, \cite{vEnRu10}, ~\cite{Opthesis09}, ~\cite{vEnKuOpRu}),  which can be subjected to Glauber dynamics, 
mixed Glauber/Kawasaki dynamics, and interacting-diffusion dynamics, 
not even necessarily Markovian (Redig, Roelly and Ruszel~\cite{ReRoRu10}), 
appliedto a large class of initial measures (e.g. Gibbs measures 
for a finite-range or an exponentially decaying interaction potential).
\item[(b)] 
Gibbs-non-Gibbs transitions can also be defined naturally for mean-field models
(see e.g. K\"ulske and Le Ny~\cite{KuLN07} for Curie-Weiss models subject to an 
independent spin-flip dynamics). In this context, much more explicit results can 
be obtained: transitions are sharp (i.e., in zero magnetic field there is a single 
time after which the measure becomes non-Gibbs and stays non-Gibbs forever, and 
in non-zero magnetic field there is a single time at which it becomes Gibbs again). 
Bad configurations can be characterized explicitly (with the interesting effect 
that non-neutral bad configurations can arise below a certain critical temperature). For further developments on mean-field results see also ~\cite{KuOp08b}, 
~\cite{ErKupr}.
\item[(c)] 
Gibbs-non-Gibbs transitions can also occur for continuous unbounded spins subject to 
independent Ornstein-Uhlenbeck processes (K\"ulske and Redig~\cite{KuRe06}), 
and for continuous bounded spins subject to independent diffusions (Van Enter 
and Ruszel~\cite{vEnRu08}, \cite{vEnRu10}), even in two dimensions where no 
static phase transitions occur.
\end{itemize}
Bad configurations can be detected by looking at a so-called \emph{two-layer 
system:} the joint distribution of the configuration at time $t=0$ and time $t>0$. 
If we condition on a particular configuration $\eta$ at time $t>0$, then the 
distribution at time $t=0$ is a Gibbs measure with an $\eta$-dependent Hamiltonian 
$H^\eta$, which is a \emph{random-field} modification of the original Hamiltonian 
$H$ of the starting measure. If, for some $\eta$, $H^\eta$ has a phase transition, 
then this $\eta$ is a bad configuration (see Fern\'andez and Pfister~\cite{FePf97}).


\subsection{Nature versus nurture}
\label{S1.2}

While these results led to a reasonably encompassing picture, we were unsatisfied 
with the strategy of the proofs for the following reason. All proofs rely on two 
fortunate facts: (1) the evolutions can be described in terms of space-time interactions; 
(2) these interactions correspond to well-studied models in equilibrium statistical 
mechanics. In particular, although the most delicate part of the analysis -- the 
proof of the onset of non-Gibbsianness -- was accomplished by adapting arguments 
developed in previous studies on renormalization transformations, the actual intuition 
that led to these results relied on entirely different arguments, based on the behavior 
of \emph{conditioned trajectories}. These intuitive arguments, already stated without proof in 
our original work~\cite{vEnFedHoRe02}, can be summarized as follows:
\begin{itemize}
\item[(I)] 
If a configuration $\eta$ is \emph{good} at time $t$ (i.e., is a point of continuity 
of the single-site conditional probabilities), then the trajectory that leads to 
$\eta$ is unique, in the sense that there is a single distribution at time $t=0$ 
that leads to $\eta$ at time $t>0$. In particular, if $t$ is small, then this 
trajectory stays close to $\eta$ during the whole time interval $[0,t]$. 
\item[(II)] 
If a configuration $\eta$ is \emph{bad} at time $t$ (i.e., is a point of essential 
discontinuity of the single-site conditional probabilities), then there are at least 
two trajectories compatible with the occurrence of $\eta$ at time $t$. Moreover, 
these trajectories can be selected by modifying the bad configuration $\eta$ 
arbitrarily far away from the origin.
\item[(III)]
Trajectories ending at a configuration $\eta$ at time $t$ are the result of a competition 
between two mechanisms:
\begin{itemize}
\item 
\emph{Nature:} The initial configuration is close to $\eta$, which is not necessarily 
typical for the initial measure, and is preserved by the dynamics up to time $t$.  
\item 
\emph{Nurture:} The initial configuration is typical for the initial measure and the 
system builds $\eta$ in a short interval prior to time $t$.
\end{itemize}
\end{itemize}

As an illustration, let us consider the low-temperature zero-field Ising model subject 
to an independent spin-flip dynamics. In~\cite{vEnFedHoRe02} we proved that the fully 
alternating configuration becomes and stays bad for large $t$. This fact can be understood 
according to the preceding paradigm in the following way.  Short times do not give the 
system occasion to perform a large number of spin-flips.  Hence, the most probable way 
to see the alternating configuration at small time $t$ is when the system started in 
a zero-magnetization-like state and the evolution kept the magnetization zero up to
time $t$. \emph{This is the nature-scenario!} For larger times $t$, a less costly 
alternative is to start in a less atypical manner, and to arrive at the alternating 
configuration following a trajectory that stays close for as long as possible to the 
unconditioned dynamical relaxation. \emph{This is the nurture-scenario!} In this 
situation, we can start either from a plus-like state or a minus-like state, as the 
difference in probabilistic cost between these two initial states is exponential in 
the size of boundary, and thus is negligible with respect to the volume cost imposed 
by a constrained dynamics. It is then possible to select between the plus-like and the 
minus-like trajectories by picking the alternating configuration in a large block, 
then picking either the all-plus or the all-minus configuration outside this block, and 
letting the block size tend to infinity.

We see that the above explanation relies on two facts:
\begin{itemize}
\item[(i)] 
The existence of a \emph{nature-versus-nurture} transition, as introduced in 
\cite{vEnFedHoRe02}.  
\item[(ii)] 
The existence of several possible trajectories (once the system is in the nurture 
regime), all starting from configurations that are typical for the initial measure 
(modulo an boundary-exponential cost). These trajectories evolve to the required bad 
configuration over a short interval prior to time $t$.
\end{itemize}


\subsection{Large deviations of trajectories}
\label{S1.3}

The goal of the present paper is to put rigor into the above qualitative suggestions. 
We propose two novel aspects: 
\begin{itemize}

\item[(1)]
 the development of a suitable large deviation theory 
for \emph{trajectories}, in order to estimate the costs of the different 
dynamical strategies; 

\item[(2)] 
the use of \emph{empirical measures} instead of \emph{configurations}, 
in order to express the conditioning at time $t$. 
\end{itemize}

For a translation-invariant spin-flip 
dynamics and a translation-invariant initial measure, nothing is lost by moving to the 
empirical measure because the bad configurations form a translation-invariant set. Instead, 
a lot is gained because, as we will show, the trajectory of the empirical measure satisfies 
a \emph{large deviation principle} under quite general conditions on the spin-flip rates 
(e.g.\ there is no restriction to high temperature). Moreover, the question of uniqueness 
versus non-uniqueness of optimal trajectories (i.e., minimizers of the large deviation 
rate function) can be posed and tackled for a large class of dynamics, which places the 
dynamical Gibbs-non-Gibbs-transition into a framework where it gains more physical 
relevance.

Here is a list of the results presented in the sequel.  

\begin{itemize}
\item[(A)] 
\emph{Existence of a large deviation principle for trajectories.} 
We apply the theory developed in Feng and Kurtz~\cite{FeKu06}, Section 8.6.  The rate 
function is the integral of the Legendre transform of the generator of the non-linear 
semigroup defined by the dynamics.  In suitably abstract terms, this generator can be 
associated to a Hamiltonian, and the rate function to the integral of a Lagrangian
(Sections~\ref{S2}--\ref{S5}).
\item[(B)]
\emph{Explicit expression for the generator of the non-linear semigroup of the dynamics.} 
These are obtained in Theorems~\ref{elprop}--\ref{elpropalt} below (Section~\ref{S3}).
\item[(C)] 
\emph{Rate functions for trajectories and associated optimal trajectories.}  
The general Legendre-transform prescription is explicitly worked out for a couple of 
simple examples, and optimal trajectories are exhibited (Sections~\ref{S4.2}--\ref{S4.3}).  
\item[(D)] 
\emph{Relation with thermodynamic potentials.}  
Relations are shown between the non-linear generator and the derivative of a ``constrained 
pressure''. Similarly, the rate function per unit time is related to the Legendre transform 
of this pressure (Section~\ref{S5.2}).
\item[(E)] 
\emph{Definition of bad measures.}  
This definition, introduced in Section~\ref{S6}, is the transcription to our more general 
framework of the notion of \emph{bad configuration} used in our original work~\cite{vEnFedHoRe02}. 
In Section~\ref{S7} we discuss the possible relations between these two notions of badness.
\end{itemize}


\subsection{Future research program}
\label{S1.4}

The results in (A)--(E) above are the preliminary steps towards a comprehensive theory of 
dynamical Gibbs-non-Gibbs transitions based on the principles outlined above. Let us conclude 
this introduction with a list of further issues which must be addressed to develop such a theory:
\begin{itemize}
\item[$\bullet$]
\emph{Definition of ``nature-trajectories'' and ``nurture-trajectories''.}
This is a delicate issue that requires full exploitation of the properties of the rate 
function for the trajectory. It must involve a suitable notion of distance between conditioned 
and unconditioned trajectories. 
\item[$\bullet$]
\emph{Relation between nature-trajectories and Gibbsianness.} 
It is intuitively clear that 
\newline
Gibbsianness is conserved for times so short that only
nature-trajectories are possible. A rigorous proof of this fact would confirm our intuition and 
would lead to alternative and less technical proofs of short-time Gibbsianness preservation.
\item[$\bullet$]
\emph{Study of nurture-trajectories.} 
We expect that nurture-trajectories start very close to unconstrained trajectories, and move
away only shortly before the end in order to satisfy the conditioning. For the case of time-reversible 
evolutions, the time it takes to get to the nurture-regime should be the same as the initial 
relaxation time to (almost) equilibrium.  
\item[$\bullet$]
\emph{Study of nature-nurture transitions.} 
Transitions from nature to nurture should happen only once for every conditioning measure (i.e.,
there should be no nature-restoration).  Natural questions are: Does the time at which these 
transitions take place depend on the conditioning measure? Is there a common time after which 
every trajectory becomes nurture?
\item[$\bullet$]
\emph{Case studies of trajectories leading to non-Gibbsianness.}  
These should determine ``forbidden regions'' in trajectory space. Natural questions are: How 
do these regions evolve? Are they monotone in time? 
\item[$\bullet$]
\emph{Relation between nurture-trajectories and non-Gibbsianness.} 
While we expect that ``all trajectories are nature'' implies Gibbsianness of the evolved measure, 
we do not expect that``some trajectories are nurture'' leads to non-Gibbsianness. Examples are 
needed to clarify this asymmetry. The case of the Ising model in non-zero field -- in which 
Gibbsianness is eventually restored -- should be particularly enlightening.
\end{itemize}


\subsection{Outline}
\label{S1.5}

Our paper is organized as follows. In Section~\ref{S2}, we consider the case of
independent spin-flips, as a warm-up for the rest of the paper. In Section~\ref{S3}, 
we compute the \emph{non-linear generator} for dependent spin-flips, which plays 
a key rol in the large deviation principle we are after. In Sections~\ref{S4}
and \ref{S5}, we compute the \emph{Legendre transform} of this non-linear 
generator, which is the object that enters into the associated rate function,
as an \emph{action integral}. In Section~\ref{S4} we do the computation for independent 
spin-flips, in Section~\ref{S5} we extend the computation to dependent spin-flips. 
In Section~\ref{S6}, we look at \emph{bad measures}, i.e., measures at time $t>0$ 
for which the optimal trajectory leading to this measure and minimizing the rate 
function is non-unique. In Section~\ref{S7}, we use these results to develop our 
large-devation view on Gibbs-non-Gibbs transtions. In Appendix~\ref{App} we 
illustrate the large deviation formalism in Feng and Kurtz~\cite{FeKu06}, which 
lies at the basis of Sections~\ref{S2}--\ref{S5}, by considering a simple example, 
namely, a Poisson random walk with small increments. This will help the reader not 
familiar with this formalism to grasp the main ideas.


\section{Independent spin-flips: trajectory of the magnetization}
\label{S2} 


\subsection{Large deviation principle}
\label{S2.1}

As a warm-up, we consider the example of Ising spins on the one-dimensional torus 
$T_N=\{1,\ldots,N\}$ subject to a rate-1 independent spin-flip dynamics. Write
$\pee_N$ to denote the law of this process. We look at the trajectory of the magnetization, 
i.e., $t \mapsto m_N(t) = N^{-1}\sum_{i=1}^N \si_i(t)$, where $\si_i(t)$ is the spin 
at site $i$ at time $t$. A spin-flip from $+1$ to $-1$ (or from $-1$ to $+1$) corresponds 
to a jump of size $-2N^{-1}$ (or $+2N^{-1}$) in the magnetization, i.e., the generator 
$L_N$ of the process $(m_N(t))_{t \geq 0}$ is given by
\be{genma}
(L_N f)(m) = \tfrac{1+m}{2}\,N\left[f\left(m-2N^{-1}\right)-f(m)\right]
+\tfrac{1-m}{2}\,N\left[f\left(m+2N^{-1}\right)-f(m)\right]
\ee
for $m\in\{-1,-1+2N^{-1},\ldots,1-2N^{-1},1\}$. If $\lim_{N\to\infty} m_N=m$ and 
$f$ is $C^1$ with bounded derivative, then 
\be{LNlim}
\lim_{N\to\infty} (L_N f)(m_N) = (Lf)(m) \quad \mbox{ with } \quad (Lf)(m) = -2mf'(m). 
\ee
This is the generator of the deterministic process $m(t) = m(0)e^{-2t}$, solving 
the equation $\dot{m}(t) = -2m(t)$ (the dot denotes the derivative with respect 
to time).

The trajectory of the magnetization satisfies a large deviation principle, i.e.,
for every time horizon $T\in (0,\infty)$ and trajectory $\gamma=(\gamma_t)_{t\in[0,T]}$,
\be{infolarge}
\pee_N\Big(\big(m_N(t)\big)_{t \in [0,T]} \approx (\gamma_t)_{t \in [0,T]}\Big) 
\approx \exp\left[- N \int_0^T L(\gamma_t,\dgam_t)\,dt\right],
\ee
where the Lagrangian $t \mapsto L(\gamma_t,\dgam_t)$ can be computed following the 
scheme of Feng and Kurtz~\cite{FeKu06}, Example 1.5. Indeed, we first compute the 
so-called \emph{non-linear generator} $H$ given by
\be{nonlingen}
(Hf)(m) = \lim_{N\to\infty} (\caH_N f)(m_N) 
\quad \mbox{ with } \quad (\caH_N f)(m_N) 
= \frac{1}{N}\,e^{-N f(m_N)}\,L_N (e^{N f})(m_N), 
\ee
where $\lim_{N\to\infty} m_N=m$. This gives
\be{nonlingen1}
(Hf)(m) = \tfrac{m+1}{2}\,(e^{-2f'(m)}-1) + \tfrac{1-m}{2}\,(e^{2f'(m)}-1),
\ee
which is of the form
\be{Hform}
(Hf)(m) = H\big(m,f'(m)\big)
\ee
with
\be{hamiltonian}
H(m,p) = \tfrac{m+1}{2}\,(e^{-2p}-1) + \tfrac{1-m}{2}\,(e^{2p}-1). 
\ee
Because $p \mapsto H(m,p)$ is convex, we have
\be{legen}
H(m,p) = \sup_{q\in\R}\,[pq-L(m,q)]
\ee
with
\be{kul}
\begin{aligned}
L(m,q) &= \sup_{p\in\R}\,[pq - H(m,p)]\\
&= \frac{q}{2} \log\left(\frac{q+\sqrt{q^2+4(1-m^2)}}{2(1+m)}\right)
-\frac12\sqrt{q^2+4(1-m^2)}+1.
\end{aligned}
\ee
Hence, using the theory developed in Feng and Kurtz~\cite{FeKu06}, Chapter 1, 
Example 1.5, we indeed have the large deviation principle in \eqref{infolarge} 
with $L(\gamma_t,\dgam_t)$ given by \eqref{kul} with $m=\gamma_t$ and $q=\dgam_t$. 


\subsection{Optimal trajectories}
\label{S2.2}

We may think of the typical trajectories $(m_N(t))_{t\in[0,T]}$ as being exponentially 
close to \emph{optimal trajectories} minimizing the \emph{action functional} $\gamma
= (\gamma_t)_{t \in [0,T]} \mapsto \int_0^T L(\gamma_t,\dgam_t)\,dt$. The optimal 
trajectories satisfy the Euler-Lagrange equations
\be{Euler}
\frac{d}{dt}\,\frac{\partial L}{\partial\dgam_t} 
= \frac{\partial L}{\partial \gamma_t}
\ee
or, equivalently, the Hamilton-Jacobi equations corresponding to the Hamiltonian
in \eqref{hamiltonian},
\be{hameq}
\dot{m} = \frac{\partial H}{\partial p}, \qquad
\dot{p} = -\frac{\partial H}{\partial m},
\ee
which gives
\be{hameq1}
\dot{m} = -m(e^{2p}+e^{-2p})+(e^{2p}-e^{-2p}), \qquad
\dot{p} = \tfrac12\,(e^{2p}- e^{-2p}).
\ee
Putting $h= \tanh(p)$ and integrating the second equation in \eqref{hameq1}, we obtain
\be{hameq2}
h(t)= C\,e^{2t}.
\ee
Using that $\mathrm{arctanh}(x)=\tfrac12\log(\tfrac{1+x}{1-x})$, we get
\be{hameq3}
\dot{m} = -m\,\frac{2+2h^2}{1-h^2}+\frac{4h}{1-h^2},
\ee
which can be integrated to yield the solution
\be{mGopt}
m(t) = C_1 e^{2t} + C_2 e^{-2t},
\ee
where the constants $C_1,C_2$ are determined by the initial magnetization
and the corresponding initial momentum. One example of an optimal trajectory 
corresponds to the dynamics \emph{starting} from an initial magnetization 
$m_0$, giving $m(t) = m_0e^{-2t}$, i.e., $C_1=0$ and $C_2=m_0$. Another 
example of an optimal trajectory is the reversed dynamics \emph{arriving} 
at magnetization $m_T$ at time $T$, giving $m(t) = m_Te^{2(t-T)}$, i.e., 
$C_2=0$ and $C_1= m_Te^{-2T}$.

Yet another example is the following. Suppose that we start the independent 
spin-flip dynamics from a measure under which the magnetization satisfies a 
large deviation principle with rate function, say, $I$, e.g.\ a Gibbs measure. 
If we want to arrive at a given magnetization $m_T$ at time $T$, then the 
optimal trajectory is given by \eqref{mGopt} with end condition $m(T)=m_T$ 
and satisfying the \emph{open-end condition} relating the Lagrangian $L$ at 
time $t=0$ to the rate function $I$ at magnetization $m=\gamma_0$ as follows:
\be{freeend}
\left[\frac{\partial L (\gamma_t,\dgam_t)}
{\partial \dot{\gamma_t}}\right]_{t=0} 
= - \left[\frac{\partial I(m)}{\partial m}\right]_{m=\gamma_0}.
\ee
This condition is obtained by minimizing $\gamma \mapsto I(\gamma_0) + \int_0^T 
L(\gamma_t,\dgam_t)\,dt$ (see Ermolaev and K\"ulske~\cite{ErKupr}).


\section{Trajectory of the empirical measure for dependent spin-flips }
\label{S3}

We will generalize the computation in Section~\ref{S2} in two directions. 
First, for independent spin-flips we are confronted with the problem that 
the rate at which the average of a local observable changes in general 
depends on the average of other observables. Second, for dependent 
spin-flips even the trajectory of the magnetization is not Markovian.
Therefore, we are obliged to consider the time evolution of all spatial 
averages jointly, i.e., the \emph{empirical measure}.


\subsection{Setting and notation}
\label{S3.1}

For $N\in\N$, let $\tend$ be the $d$-dimensional $N$-torus $(\Z/(2N+1)\Z)^d$. 
For $i,j\in\tend$, let $i+j$ denote coordinate-wise addition modulo $2N+1$. We 
consider Glauber dynamics of Ising spins located at the sites of $\tend$, i.e., 
on the configuration space $\Omega_N = \{-1,1\}^{\tend}$. We write $\Omega = 
\{-1,1\}^{\Zd}$ to denote the infinite-volume configuration space. Configurations 
are denoted by symbols like $\si$ and $\eta$. For $\si\in\Omega_N$, $\si_i$ 
denotes the value of the spin at site $i$. We write $\cM_1(\Omega)$ to denote 
the set of probability measures on $\Omega$, and similarly for $\cM_1(\Omega_N)$.

The dynamics is defined via the generator $L_N$ acting on functions $f\colon\,
\Omega_N\to\R$ as
\be{geno}
(L_N f)(\si) = \sum_{i\in\tend} c_i(\si)\,[f(\si^i)-f(\si)],
\ee
where $\si^i$ denotes the configuration obtained from $\si$ by flipping the spin
at site $i$. The rates $c_i(\si)$ are assumed to be strictly positive and 
translation invariant, i.e.,
\be{transrate}
c_i (\si) = c_0 (\tau_i \si) = c(\tau_i \si) \quad \mbox{ with } 
\quad (\tau_i\si)_j = \si_{i+j}.
\ee
We think of the dynamics with generator $L_N$ as a finite-volume version with 
periodic boundary condition of the infinite-volume generator
\be{infvol}
(Lf)(\si) = \sum_{i\in\Zd} c_i(\si)\,[f(\si^i)-f(\si)],
\ee
where now $f$ is supposed to be a local function, i.e., a function depending on a finite 
number of $\si_j$, $j\in\Zd$. We denote by $(S_t)_{t \geq 0}$ with $S_t = e^{tL}$ 
the semigroup acting on $C(\Omega)$ (the space of continuous functions on $\Omega)$) 
associated with the generator in \eqref{infvol}, and similarly $(S^N_t)_{t \geq 0}$ 
with $S^N_t= e^{tL_N}$. For $\mu\in\cM_1(\Omega)$, we denote by $\mu S_t\in\cM_1
(\Omega)$ the distribution $\mu$ evolved over time $t$, and similarly for $\mu_N 
S^N_t$ and $\mu_N\in\cM_1(\Omega_N)$.

We embed $\tend$ in $\Z^d$ by identifying it with $\Lambda_N^d=([-N,N]\cap\Z)^d$. 
Through this identification, we give meaning to expressions like $\sum_{i\in\tend} 
f(\tau_i \si)$ for $\si\in\Omega$ and $f\colon\Omega\to\R$. In this way we may also 
view local functions $f\colon\Omega\to\R$ as functions on $\Omega_N$ as soon as $N$ 
is large enough for $\Lambda_N^d$ to contain the dependence set of $f$. For a 
translation-invariant $\mu\in\cM_1(\Omega)$, we denote by $\mu_N$ its natural 
restriction to $\Omega_N$.

By the locality of the spin-flip rates, the infinite-volume dynamics is well-defined 
and is the uniform limit of the finite-volume dynamics, i.e., for every local function 
$f\colon\,\Omega\to\R$ and $t>0$,
\be{}
\lim_{N\to\infty} \|S_t^N f-S_tf\|_\infty =0.
\ee
See Liggett~\cite{Li85}, Chapters 1 and 3, for details on existence of the 
infinite-volume dynamics.


\subsection{Empirical measure}
\label{S3.2}

For $N\in\N$ and $\si\in\Omega_N$, the empirical measure associated with $\si$ is 
defined as
\be{empir}
\loc_N (\si) = \frac{1}{|\tend|} \sum_{i\in\tend} \delta_{\tau_i\si}.
\ee
This is an element of $\cM_1(\Omega_N)$ which acts on functions $f\colon\,\Omega_n
\to\R$ as
\be{}
\langle f,\loc_N\rangle = \int_{\Omega_N} f\,d\loc_N
= \frac1\norm\sum_{i\in\tend} f(\tau_i \si).
\ee
As already mentioned above, a local $f\colon\Omega\to\R$ may be considered as a 
function on $\Omega_N$ for $N$ large enough. A sequence $(\mu_N)_{N\in\N}$ 
with $\mu_N\in\cM_1(\Omega_N)$ converges weakly to some $\mu\in\cM_1(\Omega)$  
if
\be{}
\lim_{N\to\infty} \int_{\Omega_N} fd\mu_N = \int_\Omega fd\mu 
\qquad \forall\,f \mbox{ local}.
\ee
For $\si\in\Omega$, we define its periodized version $\si^N$ as $\si^N_i= \si_i$ 
for $i= (i_1,\ldots,i_d)$ with $-N\leq i_k<N$ for $k\in \{1,\ldots,d\}$, and 
$\si_i^N = \si_{i\,\mathrm{mod}\,(2N+1)}$ otherwise.

If $\mu$ is ergodic under translations, then by the locality and the translation 
invariance of the spin-flip rates also $\mu S_t$ is ergodic under translations. 
Let $\mu^N$ be the distribution of $\si^N$ when $\si$ is drawn from $\mu$. Since 
the semigroup $(S_t^N)_{t \geq 0}$ uniformly approaches the semigroup $(S_t)_{t
\geq 0}$ as $N\to\infty$, the ergodic theorem implies that
\be{conve}
\loc_N (\si^N(t)) \to \mu S_t \quad \mbox{ weakly as } \quad N\to\infty,
\ee
where $\si^N(t)$ denotes the random configuration that is obtained by evolving 
$\si^N$ over time $t$ in the process with generator $L_N$.

The deterministic trajectory $t\mapsto\mu S_t$ is the solution of the equation
\be{}
\frac{d\mu_t}{dt} = L^* \mu_t,
\ee
where $L^*$ denotes the adjoint of the generator acting on the space of finite signed 
measures $\cM(\Omega)$. Thus, we can view the convergence in \eqref{conve} as an 
infinite-dimensional law of large numbers, where the random measure-valued trajectory
$(\loc_N((\si^N(t)))_{t\in [0,T]}$ converges to the deterministic measure-valued 
trajectory $(\mu S_t)_{t\in [0,T]}$. It is therefore natural to ask for an associated 
large deviation principle, i.e., does there exist a rate function $\gamma\mapsto 
I(\gamma)$ such that
\be{}
\pee_N\Big(\big(\loc_N((\si^N(t))\big)_{t\in [0,T]} \approx \gamma\Big)
\approx \exp[-|\tend| I (\gamma)]?
\ee
Inspired by the example of the magnetization described in Section~\ref{S2}, we expect 
the answer to be yes and the rate function to be of the form
\be{intform}
I(\gamma) = \int_0^T L(\gamma_t,\dot{\gamma}_t)\,dt
\ee
for some appropriate Lagrangian $L$. In order to compute $L$, we must first find the 
generator of the non-linear semigroup.


\subsection{The generator of the non-linear semigroup}
\label{S3.3}

In our setting the non-linear generator is defined as follows:
\be{nonlin}
(H_N F)(\loc_N (\si)) = \frac1{|\tend|}\,e^{-|\tend| F(\loc_N (\si))} 
L_N\left(e^{|\tend|(F\circ\loc_N)}\right)(\si).
\ee
If the expression in \eqref{nonlin} has a limit $(HF)(\mu)$ as $N\to\infty$ when 
$\loc_N(\si)\to\mu$ weakly, then a candidate rate function can be constructed via 
Legendre transformation (see Section~\ref{S5}).

To compute the limit operator, we start with a simple function of the form
\be{}
F(\loc_N(\si)) = \langle f,\loc_N (\si)\rangle,
\ee
where $f\colon\,\Omega\to\R$ is a local function. Such $f$'s are linear combinations
of the functions
\be{haa}
H_A (\si) = \prod_{i\in A} \si_i, \qquad A\subset\Zd \mbox{ finite},
\ee
which live on $\Omega_N$ for $N$ large enough.

\bt{elprop}
For all local $f\in\Omega$ and $N$ large enough,
\be{coco}
\frac1{\norm}\,e^{-\norm\langle f,\loc_N(\si)\rangle}\,L_N 
\left(e^{|\tend|\langle f,\loc_N\rangle}\right)(\si)
= \left\langle c\left(e^{\caD_N f}-1\right),\loc_N(\si)\right\rangle,
\ee
where $\caD_N$ is the linear operator, acting on functions on $\Omega_N$, defined via
\be{dop}
\caD_N 1=0, \qquad \caD_N H_A = \sum_{r\in -A} (-2) H_{A+r} 
\mbox{ for } A\subset\tend,
\ee
where the $N$-dependence refers only to the fact that the addition $A+r$ is modulo $2N+1$.
\et

\bpr
Using the definition of the generator $L_N$ in \eqref{geno}, we write (recall 
\eqref{transrate})
\be{}
\begin{aligned}
&e^{-\norm \langle f,\loc_N (\si)\rangle}\,
L_N\left(e^{|\tend| \langle f,\loc_N\rangle}\right)(\si)\\
&\qquad = \sum_{k\in \tend} c(\tau_k \si) \left\{
\exp\left[\,\sum_{j\in \tend}\left[f\big(\tau_j (\si^k)\big)
-f\big(\tau_j(\si)\big)\right]\right]-1\right\}.
\end{aligned}
\ee
Since
\be{}
(\caD^k_Nf)(\si) = \sum_{j\in \tend} \left[f\big(\tau_j (\si^k)\big)
-f\big(\tau_j (\si)\big)\right]
\ee
is a linear operator, it suffices to prove that
\be{}
(\caD^k_N f)(\si) = (\caD_Nf)(\tau_k \si) \quad \mbox{ for } \quad  f= H_A,
\ee
where $\caD_N f$ is given by \eqref{dop} for $f=H_A$ (note that if $f=H_A$, then 
$f(\si^k)=-H_A(\si)$ for $k\in A$ and $f(\si^k)=f(\si)$ otherwise). Hence
\be{}
\begin{aligned}
(\caD^k_N H_A)(\si) 
&= \sum_{j\in\tend} \left(\prod_{i\in A} (\si^k)_{i+j} - \prod_{i\in A} \si_{i+j}\right)
= \sum_{j\in\tend} 1_{\{k-j\in A\}}\,(-2)\,\prod_{i\in A} \si_{i+j}\\
&= \sum_{j\in\tend} 1_{\{j\in -A +k\}}\,(-2)\,\prod_{i\in A} \si_{i+j}
= (-2) \sum_{r\in -A} \prod_{i\in A} \si_{i+r+k}\\ 
&= \left((-2) \sum_{r\in -A} H_{A+r}\right)(\tau_k\si).
\end{aligned}
\ee
\epr

\medskip\noindent
{\bf Remark:} Note that, in the limit as $N\to\infty$, $\caD_N$ becomes an unbounded 
operator, defined on local functions $f\colon\Omega\to\R$ via
\be{reald}
\caD 1=0, \qquad
\caD \left(\sum_{A} \alpha_A H_A\right) 
= \sum_{A} H_A \left( \sum_{r\in -A} (-2) \alpha_{A-r}\right).
\ee
The domain of $\caD$ can be extended to functions $f=\sum_A \alpha_A H_A$ for which
\be{}
\sum_{A}\sum_{r\in -A} |\alpha_{A-r}| <\infty.
\ee
The dual operator $\caD^*$ acts on $\cM(\Omega)$, the space of finite signed measures on
$\Omega$, and since $\caD 1=0$, $\caD^*$ has the measures of total mass zero as image set. The 
intuitive idea is that when the dynamics starts from the empirical measure $\mu$, after 
an infinitesimal time $t$ the empirical measure is $\mu + t \caD^*\mu + o(t)$.

\medskip\noindent
{\bf Remark:} From Theorem~\ref{elprop} it follows that, for $f=\sum_{i=1}^N \lambda_if_i$,
\be{}
L_N\left(e^{\norm\sum_{i=1}^N\lambda_i \langle f_i,\loc_N\rangle}\right)(\si)
= \norm\,e^{\norm \sum_{i=1}^n \lambda_i\langle f_i,\loc_N(\si)\rangle}
\left\langle c\left(e^{\sum_{i=1}^n \lambda_i \caD f_i}-1\right),
\loc_N(\si)\right\rangle.
\ee
The right-hand side is a function of $\loc_N$. By taking derivatives with respect to the 
variables $\lambda_i$, we see that the generator maps any function of $\loc_N$ into a 
function of $\loc_N$. This shows that $(\loc_N(\si^N(t)))_{t\geq 0}$ is a {\em Markov 
process}. Roughly speaking, this Markov process can be viewed as a random walk that 
makes jumps of size $1/\norm$ at rate $\norm$. Of course, the problem is that this 
random walk is infinite-dimensional, and therefore we cannot directly apply standard 
random-walk theory.

Theorem \ref{elprop} shows that the operator $H$ defined by
\be{}
(HF)(\mu) = \lim_{N\to\infty} (H_N F)(\loc_N(\si)) \quad \mbox{ when } \quad 
\lim_{N\to\infty} \loc_N = \mu \mbox{ weakly}
\ee
is well-defined for $F(\mu) = \langle f,\mu \rangle$. We next extend Theorem 
\ref{elprop} to $F$ of the form
\be{generf}
F(\mu) = \Psi\big(\langle f_1,\mu\rangle,\ldots,\langle f_n,\mu\rangle\big),
\ee
where $\Psi\colon\,\R^n\to\R$ is ${\mathcal C}^\infty$ with uniformly bounded derivatives 
of all orders.

\bt{elpropalt}
If $\lim_{N\to\infty} \loc_N=\mu$ and $F$ is of the form \eqref{generf}, then
(with the same notation as in \eqref{nonlin})
\be{limnonlin}
(HF)(\mu) = \lim_{N\to\infty} (H_NF)(\mu)
= \left\langle c\,\left(\exp\left[\sum_{i=1}^n \frac{\partial \Psi}{\partial x_i} 
\big(\langle f_1,\mu\rangle,\ldots,\langle f_n,\mu\rangle\big) 
\caD f_i\right]-1\right),\mu\right\rangle.
\ee
\et

\bpr
Compute
\be{bam}
\begin{aligned}
&\frac{1}{|\tend|}\,e^{-|\tend| F(\loc_N(\si))}\, L_N\,
\left( e^{|\tend| F(\loc_N)}\right)(\si)\\
&\qquad =\sum_{k\in \tend} c(\tau_k \si) \left(\exp\left[\norm\left(
F\big(\loc_N(\si^k)\big)- F\big(\loc_N (\si)\big)\right)\right]-1\right).
\end{aligned}
\ee
Next, use the fact that
\be{boem}
\langle f,\loc_N(\si^k) \rangle - \langle f,\loc_N(\si) \rangle 
= \frac{1}{\norm}\,(\caD_N f)(\tau_k (\si))
\ee
to see that
\be{baf}
\begin{aligned}
&\Psi\big(\langle f_1,\loc_N(\si^k)\rangle,\ldots,\langle f_n,\loc_N(\si^k)\rangle\big)
- \Psi\big(\langle f_1,\loc_N(\si) \rangle,\ldots,\langle f_n,\loc_N(\si)\rangle\big)\\
&\qquad = \sum_{i=1}^n \frac{\partial \Psi}{\partial x_i} 
\big(\langle f_1,\loc_N (\si)\rangle,\ldots,\langle f_n,\loc_N(\si)\rangle\big)\,
(\caD f_i)(\tau_k \si) + o\left(\frac{1}{\norm}\right).
\end{aligned}
\ee
Combine \eqref{bam} and \eqref{baf} and take the limit $N\to\infty$, to obtain 
\eqref{limnonlin}.
\epr

\medskip\noindent
{\bf Remark:} For 
\be{}
F(\mu) = \Psi\big(\langle f_1,\mu\rangle,\ldots,\langle f_n,\mu\rangle\big),
\ee
the functional derivative of $F$ with respect to $\mu$ is defined as 
\be{}
\frac{\delta F}{\delta \mu} = \sum_{i=1}^n \frac{\partial \Psi}{\partial x_i} 
\big(\langle f_1,\mu\rangle, \ldots, \langle f_n,\mu\rangle\big)\,f_i.
\ee
We may therefore rewrite \eqref{baf} as
\be{nonlinsemgen}
H(F)(\mu) = \left\langle c\,
\left(\exp\left[\caD\left(\frac{\delta F}{\delta \mu}\right)\right]-1\right),
\mu\right\rangle.
\ee


\section{The rate function for independent spin-flips}
\label{S4}

\subsection{Legendre transform}
\label{S4.1}

Having completed the computation of the non-linear generator in Section~\ref{S3},
we are ready to compute its Legendre transform. As a warm-up, we will first do 
this for independent spin-flips, i.e., when $c \equiv 1$ in \eqref{transrate}. 
In Section~\ref{S5} we will extend the calculation to general $c$, which will 
not represent a serious obstacle.

The non-linear generator in \eqref{limnonlin} is of the form
\be{legiform}
(HF)(\mu) = \caH \left(\mu,\frac{\delta F}{\delta \mu}\right),
\ee
where, for $\mu\in\cM_1(\Omega)$ and $f\colon\,\Omega\to\R$ continuous,
\be{hamilton}
\caH (\mu, f) = \left\langle c\,\left(e^{\caD f}-1\right),
\mu\right\rangle.
\ee
By the convexity of $f\mapsto \caH (\mu,f)$, we have
\be{legendre}
\caH(\mu,f) = \sup_{\alpha\in\caM(\Omega)}
\left[\int_\Omega f\,d\alpha - L(\mu,\alpha)\right]
\ee
with
\be{lagrange}
L(\mu,\alpha) = \sup_{f\in \C(\Omega)}
\left[\int_\Omega f\,d\alpha - \caH(\mu,f)\right]
\ee
the Lagrangian appearing in the large deviation rate function in \eqref{intform}. As 
explained in Feng and Kurtz~\cite{FeKu06}, Section 8.6.1.2, the representation of 
the generator in \eqref{legiform}, where $\caH(\mu,f)$ is a Legendre transform 
as in \eqref{legendre}, implies that the generator in \eqref{legiform} generates 
a non-linear semigroup, called the \emph{Nisio control semigroup}, associated with 
the function $L$ (see \cite{FeKu06}, Section 8.1).

\medskip\noindent
{\bf Remark:} The operator $\caD$ has the property
\be{}
\caD f_0 = -2f_0, \qquad f_0(\si) = \si_0,
\ee
i.e., $f_0$ is an eigenfunction of $\caD$. We recover the Hamiltonian in \eqref{hamiltonian}
(associated with the large deviation principle of the magnetization) from the 
infinite-dimensional Hamiltonian in \eqref{hamilton} via the relation
\be{}
\caH(\mu,pf_0) = H\big(\langle f_0,\mu\rangle,p\big).
\ee

\medskip\noindent
{\bf Remark:} The infinite-dimensional Hamiltonian in \eqref{hamilton} can be thought 
of as a function of the ``position'' variable $\mu$ and the ``momentum'' variable $f$.
The corresponding Hamilton-Jacobi equations read
\be{}
\dot{\mu} = \frac{\delta H}{\delta f}, \qquad
\dot{f} = -\frac{\delta H}{\delta \mu}.
\ee
These give a closed equation for $f$, because the Hamiltonian in \eqref{hamilton} is linear 
in $\mu$. If we can solve the latter equation to find $f$, then we can integrate the 
equation for $\mu$ and find the solution for $\mu$. This is precisely the same situation
-- but now infinite-dimensional -- as we encountered in \eqref{hameq1}, where the equation 
for $p$ was closed and could be integrated to give the solution for $m$. 


\subsection{Computation of the Lagrangian}
\label{S4.2}

To find $L$, the function appearing in the rate function in \eqref{intform}, we have 
to compute the Legendre transform in \eqref{lagrange}. To do so, we first consider the 
finite-dimensional analogue. We start with rates $c\equiv 1$, for which \eqref{lagrange} 
becomes
\be{findim}
\begin{aligned}
&L(\mu,\alpha) = \sup_{f\in\R^n} \left[\sum_{i=1}^n f_i \alpha_i 
- \sum_{i=1}^n \left(e^{\sum_{j=1}^n D_{ij} f_j}-1\right)\mu_i\right],\\
&\mu=(\mu_1,\ldots,\mu_n),\,\alpha=(\alpha_1,\ldots,\alpha_n),\,f=(f_1,\ldots,f_n), 
\end{aligned}
\ee
where $\mu_i\in(0,\infty)$, $\sum_{i=1}^n \mu_i=1$, $\alpha_i\in\R$, $f_i\in\R$, 
and $D_{ij}\in\R$. The matrix $D$ has the additional property that $D(\overline{1})
=\overline{0}$, where $\overline{1}$ is the vector with all components equal to $1$.
Hence $\sum_{i=1}^n (D^T\mu)_i=0$, i.e., the transposed matrix $D^T$ maps any vector 
to a vector with zero sum. For a vector $\alpha$, we say that $(D^T)^{-1}\alpha$ is 
well-defined if there exists a unique vector $\nu=\nu(\alpha)$ with sum equal to $1$ 
such that $D^T \nu =\alpha$. For two column vectors $\alpha,\beta\in \R^n$, let 
$\alpha\beta$ be the vector with components $\alpha_i\beta_i$, $\alpha/\beta$ the 
vector with components $\alpha_i/\beta_i$. For $f\colon\,\R\to\R$, write $f(\alpha)$ 
to denote the vector with components $f(\alpha_i)$. Then the equation for the maximizer 
$f=f^*$ of \eqref{findim} becomes
\be{findimeq}
\alpha_k = \sum_{i=1}^n \mu_i\,e^{\sum_{j=1}^n D_{ij} f^*_j}\,D_{ik},
\quad k=1,\ldots,n,
\ee
which in vector notation reads
\be{}
\alpha= D^T(\mu e^{Df^*}).
\ee
If $(D^T)^{-1}\alpha$ is well-defined, then we can rewrite the latter equation as
\be{}
Df^* = \log\left[\frac{(D^T)^{-1}\alpha}{\mu}\right],
\ee
and for this $f^*$ we have
\be{var1}
\sum_{i=1}^n f^*_i \alpha_i = \langle f^*, \alpha\rangle
= \left\langle \log\left[\frac{(D^T)^{-1}\alpha}{\mu}\right],
(D^T)^{-1}\alpha\right \rangle
\ee
and
\be{var2}
\sum_{i=1}^n \left( e^{\sum_{j=1}^n D_{ij} f^*_j}-1\right)\mu_i=0,
\ee
because the total mass of $\mu$ and $(D^T)^{-1}$ are both equal to $1$. Hence, 
inserting \eqref{var1} and \eqref{var2} into \eqref{findim}, we 
obtain the expression
\be{lagrafindim}
L(\mu,\alpha) = \left\langle \log\left[\frac{(D^T)^{-1}\alpha}{\mu}\right],
(D^T)^{-1}\alpha\right\rangle,
\ee
which is the \emph{relative entropy} of $(D^T)^{-1}\alpha$ with respect to $\mu$. 
The intuition behind \eqref{lagrafindim} is that $L(\mu,\alpha)$ is the cost under 
the Markovian evolution for the initial measure to have derivative $\alpha$ at time 
zero.

Let us next consider the infinite-dimensional version of the above computation. First, 
for $\alpha\in\cM(\Omega)$ with total mass zero, we declare $(\caD^*)^{-1}\alpha=\nu$ 
to be well-defined if there exists a probability measure $\nu$ such that, for
all $f$ in the domain of $\caD$,
\be{}
\langle\nu,\caD f \rangle = \langle\alpha,f\rangle.
\ee
If $\alpha$ is translation-invariant, then also $(\caD^*)^{-1}\alpha$ is 
translation-invariant. For translation-invariant $\mu,\nu\in\cM(\Omega)$,
we denote by $s(\nu|\mu)$ the relative entropy density of $\nu$ with respect to $\mu$, 
i.e.,
\be{raaf}
s(\nu|\mu) = \lim_{N\to\infty} \frac{1}{|\tend|} \sum_{\si_{\tend}} \nu(\si_{\tend}) 
\log\left[\frac{\nu(\si_{\tend})}{\mu(\si_{\tend})}\right].
\ee
Note that this limit does not necessarily exist. But if $\mu$ is a Gibbs measure, then 
for all translation-invariant $\nu$ both $s(\nu|\mu)$ and $s(\nu|\mu_t)$ exist, where
$\mu_t$ is $\mu$ evolved over time $t$ (see van Enter, Fern\'andez and Sokal~\cite{vEnFeSo93}, 
Le Ny and Redig~\cite{LNRe04}). The rate function which is the analogue of \eqref{lagrafindim} 
is now given by
\be{lagra}
L(\mu,\alpha)= s\left((\caD^*)^{-1}\alpha|\mu\right)
\ee
with the same interpretation as for \eqref{lagrafindim}: $(\caD^*)^{-1}\alpha$ produces 
derivative $\alpha$ at time zero for the trajectory of the empirical measure, and its cost 
is the relative entropy density of this measure with respect to the initial measure $\mu$.


\subsection{Optimal trajectories}
\label{S4.3}

In order to gain some intuition about the rate function corresponding to the Lagrangian 
in \eqref{lagra}, we identify two easy optimal trajectories. 

First, we consider a trajectory that starts from a product measure $\nu_{x_0}$ and ends 
at a product measure $\nu_{x_t}$ with $x_t= x_0 e^{-2t}$. The typical trajectory 
is then simply the product-measure-valued trajectory $\gamma_t=\nu_{x_t}$ with 
$x_t= x_0 e^{-2t}$. We can easily verify that this trajectory has zero cost. Indeed,
$\langle\gamma_t,H_A\rangle = x_t^{|A|}$, and hence $\langle\dot{\gamma_t},H_A\rangle 
= |A| x_t^{|A|-1} \dot{x}_t$. On the other hand, $\langle\caD^*(\gamma_t),H_A\rangle 
= -2|A| x_t^{|A|}$ and, since $\dot{x}_s= -2 x_t$, we thus see that $\langle\dot{\gamma_t},
H_A\rangle = \langle\caD^* (\gamma_t),H_A\rangle$. Therefore $(\caD^*)^{-1}(\dot{\gamma}_t)
= \gamma_t$, and \eqref{lagra} gives
\be{}
L(\gamma_t,\dot{\gamma}_t) = s\left((\caD^*)^{-1}(\dot{\gamma}_t)|\gamma_t\right)
= s(\gamma_t|\gamma_t) = 0. 
\ee
Note that this is the only product-measure-valued trajectory that has zero cost.
Indeed, if $\gamma_t=\nu_{x_t}$ has zero cost, then from the requirement that 
$\langle\dot{\gamma_t},H_A\rangle = \langle\caD^*(\gamma_t),H_A\rangle
= -2|A| x_t^{|A|}$ we find that $\dot{x}_t =-2 x_t$. For a general starting measure 
$\mu$, the trajectory that has zero cost satisfies $\langle\dot{\gamma}_t,H_A\rangle 
= -2|A|\,\langle\gamma_t,H_A\rangle$, which has as solution $\langle\gamma_t,H_A\rangle 
= \langle\mu,H_A\rangle\,e^{-2|A|t}$, corresponding to the Markovian independent spin-flip 
evolution started from $\mu$. Note that, for a general trajectory $\gamma$, 
$\langle(\caD^*)^{-1}(\dot{\gamma}_t),H_A\rangle = -2|A|\langle\dot{\gamma}_t,H_A\rangle$.

Second, we consider the case where $\mu=\mu_y$ is a product measure with
\be{}
\langle\mu_y,H_A\rangle = y^{|A|}, \qquad -1<y<1,
\ee
and $\alpha = \alpha_x$ is the derivative at time zero of another product measure,
i.e.,
\be{}
\langle\alpha_x,H_A\rangle = -2|A|x^{|A|}, \qquad -1<x<1.
\ee
In that case $D^*\alpha = \nu_x$ with $\nu_x$ the translation-invariant product measure 
with $\langle\nu_x,H_A\rangle= x^{|A|}$. The latter follows from the identity
\be{}
\left\langle\left[\sum_{i\in A} H_A (\si^i)- H_A(\si)\right],\nu_x\right\rangle 
= -2|A|x^{|A|},
\ee
and the rate function becomes
\be{}
L(\mu_y,\alpha_x) = \frac{1+x}{2} \log\left(\frac{1+x}{1+y}\right)
+ \frac{1-x}{2} \log\left(\frac{1-x}{1-y}\right).
\ee


\section{The rate function for dependent spin-flips}
\label{S5}


\subsection{Computation of the Lagrangian}
\label{S5.1}

For general spin-flip rates $c$ in \eqref{transrate}, let us return to the matrix calculation 
in \eqref{findim} and \eqref{findimeq}. Equation \eqref{findim} has to be replaced by
\be{findimc}
L(\mu,\alpha) = \sup_{f\in\R^n} \left[\sum_{i=1}^n f_i \alpha_i 
- \sum_{i=1}^n c_i \left(e^{\sum_{j=1}^n D_{ij} f_j}-1\right)\mu_i\right],
\ee
where $c_i>0$, $i=1,\ldots,n$. Put $C_\mu=\sum_{i=1}^n c_i\mu_i$. In the calculation with 
$c_i=1$, $i=1,\ldots,n$, this ``total mass'' does not depend on $\mu$ and is equal to 1. 
Now, however, it becomes a normalization that depends on $\mu$. We say that $(D^T)^{-1}
(\alpha,\mu)$ is well-defined if there exists a non-negative vector $\nu=\nu(\alpha,\mu)
=(\nu_1\ldots,\nu_n)$ with sum $C_\mu$ such that $D^T\nu = \alpha$. The analogue of 
\eqref{lagrafindim} reads
\be{lagrafindimc}
L(\mu,\alpha) = \left\langle\log\left[\frac{(D^T)^{-1}(\alpha,\mu)}{\mu}\right],
(D^T)^{-1}(\alpha,\mu)\right\rangle.
\ee

In order to find the analogue of this expression in the infinite-dimensional setting,
we proceed as follows. For two finite positive measures $\mu,\nu$ of equal total 
mass $M$, we define $S(\mu|\nu)$ to be the relative entropy density of the probability 
measures $\mu/M, \nu/M$, i.e., $S(\mu|\nu)=s(\nu/M|\mu/M)$. For $\mu\in\cM_1(\Omega)$, 
we define the $c$-modification of $\mu$ as the positive measure defined via $\int_\Omega 
f(\si)d\mu_c(\si)=\int_\Omega c(\si)f(\si)d\mu(\si)$. For a signed measure of total mass 
zero and $\mu\in\cM_1(\Omega)$, we say that $(\caD^*)^{-1}(\alpha,\mu,c)=\nu$ is 
well-defined if there exists a positive measure $\nu$ of total mass equal to that of 
$\mu_c$ such that $\caD^*(\nu)=\alpha$. Then the analogue of \eqref{lagrafindimc} becomes
\be{lagrainfdimc}
L(\mu,\alpha)= s\left((\caD^*)^{-1}(\alpha,\mu,c)|\mu_c\right).
\ee


\subsection{The non-linear semigroup and its relation with relative entropy}
\label{S5.2}

The non-linear semigroup with generator \eqref{nonlin} is defined as follows. Let 
$\caP^{\mathrm{inv}}(\Omega)$ be the set of translation-invariant probability measures 
on $\Omega$. For local functions $f_1,\ldots,f_n\colon\,\Omega\to\R$ and a 
$\ce^\infty$-function $\Psi\colon\,\R^n\to\R$, we define an associated function 
$F^{f_1,\ldots,f_N}_\Psi\colon\,\caP^{\mathrm{inv}}(\Omega)\to\R$ via 
\be{}
F^{f_1,\ldots,f_n}_\Psi (\mu) 
= \Psi\left(\int_\Omega f_1 d\mu,\ldots,\int_\Omega f_n d\mu\right).
\ee
Since $\langle f_i,\loc_N\rangle$ is well-defined for $N$ large enough, we can 
define $F^{f_1,\ldots,f_n}_\Psi (\loc_N)$ for $N$ large enough as well. This allows
us to define a non-linear semigroup $(V(t))_{t\geq 0}$ via
\be{semigrogi}
\big(V(t) F^{f_1,\ldots,f_n}_\Psi\big)(\mu) = \lim_{N\to\infty} \frac{1}{\norm}\, 
\log\E_{\sigma^N}\left(\exp\left[\norm F^{f_1,\ldots,f_n}_\Psi\big(\loc_N(\si^N(t))\big)
\right]\right), 
\ee
where $\E_{\sigma^N}$ denotes expectation with respect to the law of the process 
starting from $\sigma^N$, and the limit is taken along a sequence of configurations 
$(\si^N)_{N\in\N}$ with $\si^N\in\Omega_N$ such that the associated empirical measure 
$\loc_N(\si^N)$ converges weakly to $\mu$ as $N\to\infty$. If $V(t)$ exists, then it 
defines a non-linear semigroup, and the generator of $V(t)$ is given by 
\eqref{nonlinsemgen}.

Conversely, if $H$ in \eqref{nonlinsemgen} generates a semigroup, then this must be
$(V(t))_{t\geq 0}$. The fact that this semigroup is well-defined is sufficient 
to imply the large deviation principle for the trajectory of the empirical measure
(Feng and Kurtz~\cite{FeKu06}, Theorem 5.15). Technically, the difficulty consists 
in showing that the generator in \eqref{nonlinsemgen} actually generates a semigroup.

We now make the link between the non-linear semigroup, its generator and some 
familiar objects of statistical mechanics, such as pressure and relative entropy 
density. 

\bd{cpres}
The constrained pressure at time $t$ associated with a function $f\colon\,\Omega\to\R$ 
and a Gibbs measure $\mu\in\caP^{\mathrm{inv}}(\Omega)$ is defined as
\be{pres}
p_t(f|\mu) = \lim_{N\to\infty} \frac{1}{\norm}\,
\log\E_{\si^N}\left(e^{\sum_{x\in \T_N} \tau_x f(\si_t)}\right),
\ee
where the limit is taken along a sequence of configurations $(\si^N)_{N\in\N}$
with $\si^N\in\Omega_N$ such that the associated empirical measure $\loc_N(\si^N)$ 
converges weakly to $\mu$ as $N\to\infty$.
\ed
In particular, $p_0 (f|\mu)= \int_\Omega f d\mu$. The relation between the non-linear 
semigroup in \eqref{semigrogi} and the constrained pressure at time $t$ reads
\be{pressem}
\big(V(t)\langle f,\cdot\rangle\big)(\mu) = p_t(f|\mu).
\ee

The pressure at time $t$ is defined as
\be{presst}
p(f|\mu_t) = \lim_{N\to\infty} \frac{1}{\norm}\,
\log\E_{\mu}\left(e^{\sum_{x\in \T_N} \tau_x f(\si_t)}\right).
\ee
This is well-defined as soon as the dynamics starts from a Gibbs measure $\mu_0=\mu$ (see
Le Ny and Redig~\cite{LNRe04}). The relation between the pressure and the constrained 
pressure reads
\be{pressrel}
p(f|\mu_t) = \sup_{\nu\in\caP^{\mathrm{inv}}(\Omega)}\,[p_t(f|\nu) -s(\nu|\mu)].
\ee
On the other hand, the pressure at time $t$ is the Legendre transform of the relative
entropy density with respect to $\mu_t$, i.e., 
\be{legentrop}
p(f|\mu_t) = \sup_{\nu\in\caP^{\mathrm{inv}}(\Omega)}\,
\left[\int_\Omega f d\nu- s(\nu|\mu_t)\right],
\ee
where the relative entropy density $s(\nu|\mu_t)$ exists because $\mu_t$ is 
\emph{asymptotically decoupled} (see Pfister~\cite{Pf??}) as soon as 
$\mu_0=\mu$ is a Gibbs measure (see Le Ny and Redig~\cite{LNRe04}).

The relation between the non-linear generator and the constrained pressure is now 
as follows. Define the Legendre transform of the constrained pressure as
\be{legendrepress}
p^*_t (\nu|\mu) = \sup_{f\in\C(\Omega)} \left[\int f d\nu - p_t (f|\mu)\right].
\ee
Then the relation with the Hamiltonian in \eqref{hamilton} and the Lagrangian in
\eqref{lagrainfdimc} is 
\be{hampress}
H(\mu,f) = \left[\frac{d}{dt}\,p_t(f|\mu)\right]_{t=0}
\ee
and
\be{lagpress}
L(\mu,\alpha) = \lim_{t\to 0} \frac{1}{t}\,p^*_t(\mu+t\alpha|\mu).
\ee

\medskip\noindent
{\bf Remark:} The operator $\caD$, acting on the space $C(\Omega)$ of continuous functions
on $\Omega$, has a dual operator $\caD^*$, acting on the space $\cM(\Omega)$ of finite signed 
measures on $\Omega$, defined via
\be{}
\langle f,\caD^*\mu\rangle = \langle\caD f,\mu\rangle.
\ee
In order to gain some understanding for $\caD^*$ (which will be useful later on), 
we first compute $\caD^*$ for a Gibbs measure $\mu\in\caP^{\mathrm{inv}}(\Omega)$.
Without loss of generality we may assume that the interaction potential of $\mu$ 
is a sum of terms of the form $\Phi(A,\si)=J_A H(A,\si)$, $A\subset\Zd$ finite, where $J_A$ is 
translation invariant, i.e., $J_{A+k} = J_A$, $k\in\Zd$. We also assume absolute 
summability, i.e.,  
\be{}
\sum_{A\ni 0} |J_{A}| <\infty.
\ee
Remember that
\be{}
(\caD f)(\si) = \sum_{j\in\Zd} \left[f(\tau_j(\si^0))- f(\tau_j(\si))\right].
\ee
Therefore, for the Gibbs measure $\mu$ under consideration, we have
\be{gibbsd}
\langle\mu,\caD f\rangle = \int_\Omega\left(\sum_{j\in\Zd}\frac{d\mu^0}{d\mu}
\circ\tau_{-j}-1\right)\,fd\mu,
\ee
where $\mu^0$ denotes the distribution of $\si^0$ when $\si$ is distributed according 
to $\mu$. Note that the sum in the right-hand side of \eqref{gibbsd} is formal, i.e., 
the integral is well-defined due to the multiplication with the local function $f$. 
In terms of $J_A$, $A\subset\Zd$ finite, we have
\be{}
\left(\sum_{j\in\Zd} \frac{d\mu^0}{d\mu}\circ\tau_{-j}-1\right)(\si)
= \sum_{j\in\Zd} \left(e^{-\sum_{A\ni 0}-2J_A H(A-j,\si)}-1\right),
\ee
where, once again, this expression is well-defined only after multiplication with a 
local function and integrated over $\mu$.


\section{Bad empirical measures}
\label{S6}

In Section \ref{S7} we will see what consequences the large deviation principle for 
the trajectory of the empirical measure derived in Sections \ref{S3} and \ref{S5} 
has for the question of Gibbs versus non-Gibbs. This needs the notion of bad 
empirical measure, which we define next.
 
If we start our spin-flip dynamics from a Gibbs measure $\mu\in\caP^{\mathrm{inv}}
(\Omega)$, then a probability-measure-valued trajectory $\gamma=(\gamma_t)_{t\in[0,T]}$ 
has cost
\be{}
\caI_\mu (\gamma) = s(\gamma_0|\mu) + \int_0^T L(\gamma_t,\dot{\gamma}_t)\,dt,
\ee
where the term $s(\gamma_0|\mu)$ is the cost of the initial distribution $\gamma_0$. 
We are interested in the minimizers of $\caI_\mu (\gamma)$ over the set of trajectories 
$\gamma$ that {\em end} at a given measure $\nu$. Let
\be{}
\begin{aligned}
K_T (\mu',\nu) &= \inf_{\gamma\colon\,\gamma_0 = \mu',\gamma_T=\nu}
\int_0^T L(\gamma_t,\dot{\gamma}_t)\,dt\\
&= - \lim_{N\to\infty} \frac{1}{\norm}\,
\log\pee_\mu\left(\loc_N(\si_T)=\nu|\loc_N(\si_0)=\mu'\right).
\end{aligned}
\ee
Then $e^{\norm -K_T(\mu',\nu)}$ can be thought of as the transition probability
for the empirical measure $\loc_N$ to go from $\mu'$ to $\nu$, up to factors of order 
$e^{o(\norm)}$. Hence
\be{sK}
\begin{aligned}
&-\frac{1}{\norm}\,\log\pee_\mu\left(\loc_N(\si_0)= \mu'|\loc_N (\si_T)=\nu\right)\\
&\qquad = [s(\mu'|\mu) + K_T (\mu',\nu)] 
- \inf_{\mu'\in\caP^{\mathrm{inv}}(\Omega)}\,[s(\mu'|\mu) + K_T (\mu',\nu)].
\end{aligned}
\ee
Let $M^*(\mu,\nu)$ be the set of probability measure $\mu'$ for which the infimum 
in the right-hand side of \eqref{sK} is attained. We can then think of each element 
in this set as a typical empirical measure at time $t=0$ given that the empirical 
measure at time $T$ is $\nu$. When $M^*$ is a singleton, we denote its unique element 
by $\mu^*(\mu,\nu)$.

\bd{bad}
(a) A measure $\nu$ is called bad at time $t$ if $M^*(\mu,\nu)$ contains at least 
two elements $\mu_1$ and $\mu_2$ and there exist two sequences $(\nu^1_n)_{n\in\N}$
and $(\nu^2_n)_{n\in\N}$, both converging to $\nu$ as $n\to\infty$, such that
$\mu^*(\mu,\nu^1_n)$ converges to $\mu_1$ and $\mu^*(\mu,\nu^2_n)$ converges
to $\mu_2$.\\ 
(b) A measure $\nu$ that is bad at time $t$ has at least two possible histories,
stated as a two-layer property: seeing the measure $\nu$ at time $t$ is compatible 
(in the sense of optimal trajectories) with two different measures at time $t=0$.
\ed

\noindent
Badness of a measure can be detected as follows.

\bp{badprop}
A measure $\nu$ is bad at time $t$ if there exists a local function $f\colon\Omega\to\R$, 
two sequences $(\nu^1_n)_{n\in\N}$ and $(\nu^2_n)_{n\in\N}$ both converging to
$\nu$, and an $\epsi>0$ such that
\be{baddefi}
\Big|\E\Big(f(\si(0))~|~\loc_N(\si(t))=(\nu^1_n)_N\Big)
-\E\Big(f(\si(0))~|~\loc_N(\si(t))=(\nu^2_n)_N\Big)\Big| > \epsi \quad \forall\,N,n\in\N,
\ee
where $(\nu_n)_N$ denotes the projection of $\nu_n$ on $\tend$.
\ep


\section{A large deviation view on dynamical Gibbs-non-Gibbs transitions}
\label{S7}

In van Enter, Fern\'andez, den Hollander and Redig~\cite{vEnFedHoRe02} we studied 
the evolution of a Gibbs measure $\mu$ under a high-temperature spin-flip dynamics. 
We showed that the Gibbsianness of the measure $\mu_t$ at time $t>0$ is equivalent 
to the absence of a phase transition in the double-layer system. More precisely, 
conditioned on end configuration $\eta$ at time $t$, the distribution at time $t=0$ 
is a Gibbs measure $\mu^\eta$ with $\eta$-dependent formal Hamiltonian
\be{}
 H^\eta_t (\si, \eta) = H(\si) + h_t \sum_{i\in\Zd} \si_i \eta_i,
\ee
where $t\mapsto h_t$ is a monotone function with $\lim_{t\downarrow 0} h_t=\infty$ 
and $\lim_{t\to\infty} h_t=0$. If the double-layer system has a phase transition for 
an end configuration $\eta$, then $\eta$ is called bad. In that case $\eta$ is an 
essential point of discontinuity for any version of the conditional probability 
$\mu_t (\si_\la=\cdot\,|\si_{\la^c})$, $\la\subset\Zd$ finite.

The relation between the double-layer system and the trajectory of the empirical 
distribution is as follows. Suppose that the double-layer system has no phase 
transition for any end configuration $\eta$. If we condition the empirical 
measure at time $t>0$ to be $\nu$, then (by further conditioning on the configuration 
$\eta$ at time $t>0$) we conclude that at time $t=0$ we have the measure $\int_\Omega
\mu^\eta\nu(d\eta)$. Hence the optimizing trajectory is unique. Conversely, if there exist 
a bad configuration $\eta$, then (because of the translation invariance of the 
initial measure and of the dynamics) all translates of $\eta$ are bad also. Hence
we expect that a translation-invariant measure $\nu$ arising as any weak limit 
point of $\norm^{-1}\sum_{x\in\tend}\delta_{\tau_x\eta}$ is bad also.

As an example, let us consider a situation studied in \cite{vEnFedHoRe02}. The
dynamics starts from $\mu^{+}_\beta$, the low-temperature plus-phase of the Ising 
model with zero magnetic field, and evolves according to independent spin-flips.
Then, from some time onwards, the alternating configuration $\eta_{\mathrm{alt}}
(x)=(-1)^{\sum_{i=1}^d |x_i|}$ becomes bad. The same is true for $-\eta_{\mathrm{alt}}$,
and so the translation-invariant measure 
\be{nusplit}
\nu = \tfrac12\left(\delta_{\eta_{\mathrm{alt}}}+\delta_{-\eta_{\mathrm{alt}}}\right)
\ee
has the property that, for $\nu$-a.e.\ configuration $\eta$, the double-layer system 
has a phase transition when the end configuration is $\eta$. Moreover, the Hamiltonian 
$H^\eta_t$ has a plus-phase $\mu^+_\eta$ and a minus-phase $\mu^-_\eta$. Therefore, 
when we condition on the empirical measure in \eqref{nusplit} we get two possible optimal 
trajectories, one starting at $\frac12(\mu^+_\eta + \mu^{+}_{-\eta})$ and one starting 
at $\frac12(\mu^-_\eta + \mu^{-}_{-\eta})$. To realize the approximating measures of 
Proposition~\ref{badprop}, we choose $\nu^1_n,\nu^2_n$ to be the randomized versions of 
$\nu$ where we first choose a configuration according to $\nu$ and then 
independently flip spins with probability $1/n$, 
to change either from minus to plus or  stay plus if it was plus to begin with, 
respectively to change to minus or stay minus.
Clearly, by the FKG-inequality, when conditioning 
on $\nu^1_n$, respectively, $\nu^2_n$ as empirical distribution, we get a measure at 
time $t=0$ that is above $\mu^+_\eta + \mu^+_{-\eta}$, respectively, below $\mu^-_\eta+
\mu^-_{-\eta}$. Hence \eqref{baddefi} holds with $f(\si)=\si_0$, and $\nu$ is bad.


\appendix

\section{A simple example of the Feng-Kurtz formalism}
\label{App}
 
\subsection{Poisson walk with small increments}
\label{A1}

In order to introduce the general formalism developed in Feng and Kurtz~\cite{FeKu06},
let us consider a simple example where computations are simple yet the fundamental 
objects of the general theory already appear naturally.

Let $X_N=(X_N(t))_{t \geq 0}$ be the continuous-time random walk on $\R$ that jumps 
$N^{-1}$ forward at rate $bN$ and $-N^{-1}$ backward at rate $dN$, with $b,d \in 
(0,\infty)$. This is the Markov process with generator
\begin{equation}
\label{poissongen}
(L_N f)(x) = bN \left[f\left(x+N^{-1}\right) - f(x)\right] 
+ dN \left[f\left(x-N^{-1}\right) - f(x)\right]. 
\end{equation}
Clearly, if $\lim_{N\to\infty} X_N(0)=x\in\R$, then 
\begin{equation}
\lim_{N\to\infty} X_N(t) = x + (b-d)t, \qquad t>0,
\end{equation}
i.e., in the limit as $N\to\infty$ the random process $X_N$ becomes a deterministic 
process $(x(t))_{t\geq 0}$ solving the limiting equation
\begin{equation}
\dot{x} =(b-d), \qquad x(0)=x.
\end{equation}

For all $N\in \N$, we have
\begin{equation}
X_N(t) = N^{-1}\,\left[\cN^+(Nbt)-\cN^-(Ndt)\right] = \sum_{i=1}^N (X_i^{bt}-Y_i^{dt})
\end{equation}
with $\cN^+=(\cN^+(t))_{t\geq 0}$ and $\cN^-=(\cN^-(t))_{t\geq 0}$ independent rate-1 
Poisson processes, and $X_i^t$, $Y_i^t$, $i=1,\dots,N$, independent Poisson random 
variables with mean $bt$, respectively, $dt$. Consequently, we can use Cram\'er's
theorem for sums of i.i.d.\ random variables to compute
\begin{equation}
\label{poissonrate}
I(at) = \lim_{N\to\infty} \frac{1}{N} \log\pee_N
\big(X_N (t) = at \mid  X_N(0)=0\big)
= \sup_{\lambda\in\R} \big[at\lambda - F(\lambda)\big],
\end{equation}
where
\begin{equation}
F(\lambda) = \lim_{N\to\infty} \frac{1}{N} \log \E_N\left(e^{\lambda N X_N(t)}\right) 
= b \big(e^{\lambda}-1\big) + d \big(e^{-\lambda}-1\big).
\end{equation}
Thus, we see that
\begin{equation}
I(at) = t L(a)
\end{equation}
with
\begin{equation}
\label{poissonlagra}
L(a) = \sup_{\lambda \in\R} \left[a\lambda-b\big(e^{\lambda}-1\big) 
- d \big(e^{-\lambda}-1\big)\right].
\end{equation}
Using the property that the increments of the Poisson process are independent over disjoint time intervals,
we can now compute
\begin{equation}
\label{ratecomputation}
\begin{aligned}
&\lim_{N\to\infty} \frac{1}{N}\log\pee_N\Big((X_N (t))_{t\in [0,T]} 
\approx (\gamma_t)_{t\in [0,T]}\Big)\\
&\qquad =
\lim_{n\to\infty}
\sum_{i=1}^n \lim_{N\to\infty} \frac{1}{N}\log\pee_N\Big(X_N(t_i)-X_N(t_{i-1}) 
\approx \dot{\gamma}_{t_{i-1}}(t_i-t_{i-1})\Big)\\
&\qquad = \lim_{n\to\infty} \sum_{i=1}^n (t_i-t_{i-1})\,L\big(\dot{\gamma}_{t_{i-1}}\big)
= \int_0^T  L\big(\dot{\gamma}_t\big)\,dt,
\end{aligned}
\end{equation}
where $L$ is given by \eqref{poissonlagra} and $t_i$, $i=1,\ldots,n$, is a partition of the 
time interval $[0,T]$ that becomes dense in the limit as $n\to\infty$.

We see from the above elementary computation that, in the limit as $N\to\infty$,
\begin{equation}
\label{deville}
\pee_N\Big((X_N(t))_{t\in[0,T]} \approx (\gamma(t))_{t\in[0,T]}\Big) 
\approx \exp\left[-N\int_0^T L(\gamma_t,\dot{\gamma}_t)\,dt\right],
\end{equation}
where the \emph{Lagrangian} $L$ only depends on the second variable, namely,
\begin{equation}
L(\gamma_t,\dot{\gamma}_t)= L(\dot{\gamma}_t)
\end{equation}
with $L$ given by \eqref{poissonlagra}. We interpret \eqref{deville} as follows: if the 
trajectory is not differentiable at some time $t \in [0,T]$, then the probability in the 
left-hand side of \eqref{deville} decays superexponentially fast with $N$, i.e., 
\begin{equation}
\lim_{N\to\infty} \frac{1}{N}\log\pee_N\Big((X_N(t))_{t \in [0,T]} 
\approx (\gamma_t)_{t \in [0,T]}\Big) = -\infty, 
\end{equation}
and otherwise it is given by the formula in \eqref{deville} (read in the standard
large-deviation language).

The \emph{Lagrangian} in \eqref{poissonlagra} is the Legendre transform of the \emph{Hamiltonian}
\begin{equation}
H(\lambda) = b\big(e^{\lambda}-1\big) - d \big(e^{-\lambda}-1\big).
\end{equation}
This Hamiltonian can be obtained from the generator in \eqref{poissongen} as follows:
\begin{equation}
H(\lambda) = \lim_{N\to\infty} \frac{1}{N}\,e^{-Nf_\lambda(x)}\,
\big(L_N e^{N f_\lambda}\big)(x), \qquad f_\lambda(x) = \lambda x.
\end{equation}
More generally, by considering the operator
\begin{equation}
\label{poissonnonlingen}
(\caH f)(x) = \lim_{N\to\infty} \frac{1}{N}\,e^{-N f(x)}\,\left(L_N e^{Nf}\right)(x) 
= b \big(e^{f'(x)}-1\big) - d \big(e^{-f'(x)}-1\big),
\end{equation}
we see that the Hamiltonian equals
\begin{equation}
H(\lambda) = (\caH f_\lambda)(x),
\end{equation}
and that, by the convexity of $\lambda \mapsto H(\lambda)$,
\begin{equation}
(\caH f)(x) = H(f'(x)) = \sup_{a\in \R} [a f'(x)- L(a)].
\end{equation}
The operator $\caH$ is called the \emph{generator of the non-linear semigroup}.


\subsection{The scheme of Feng and Kurtz}
\label{A2}

The scheme that produces the Lagrangian in \eqref{poissonlagra} from the operator in 
\eqref{poissonnonlingen} actually works in much greater generality. Consider a sequence 
of Markov processes $X=(X_N)_{N\in\N}$ with $X_N=(X_N(t))_{t\geq 0}$, living on a common 
state space (like $\R$, $\R^d$ or a space of probability measures). Suppose that $X_N$ 
has generator $L_N$ and in the limit as $N\to\infty$ converges to a process $(x(t))_{t\geq 0}$, 
which can be either deterministic (as in the previous example) or stochastic. We want to 
identify the Lagrangian controlling the large deviations of the trajectories:
\begin{equation}
\label{pathspacelargedevpri}
\pee_N\Big((X_N (t))_{t \in [0,T]} \approx (\gamma_t)_{t \in [0,T]}\Big) 
\approx \exp\left[-N \int_0^T L(\gamma_t,\dot{\gamma}_t)\,dt\right].
\end{equation}
Omitting technical conditions, we see that this can be done in four steps:
\begin{enumerate}
\item 
Compute the generator of the non-linear semigroup 
\begin{equation}
\label{nonlingengen}
(\caH f)(x) = \lim_{N\to\infty} \frac{1}{N}\,e^{-N f(x)}\,\left(L_N e^{N f}\right)(x).
\end{equation}
\item 
Look for a function $H(x,p)$ of two variables such that 
\begin{equation}
(\caH f)(x) = H(x,\nabla f(x)).
\end{equation}
What $\nabla f$ means depends on the context: on $\R^d$ it simply is the gradient of $f$, 
while on an infinite-dimensional state space it is a functional derivative.
\item 
Express the function $H$ as a Legendre transform:
\begin{equation}
H(x,p) = \sup_{p} \left[\langle p,\lambda\rangle - L(x,\lambda)\right].
\end{equation}
What $\langle\cdot\rangle$ means also depends on the context: on $\R^d$ it simply is
the inner product, while in general it is a natural pairing between a space and its 
dual, such as $\langle f,\mu\rangle = \int f d\mu$.
\item 
The Lagrangian in \eqref{pathspacelargedevpri} is the function $L$ with $x=\gamma_t$ 
and $\lambda =\dot{\gamma}_t$.
\end{enumerate}


\end{document}